\documentclass[12pt]{article}
\usepackage{amssymb}
\usepackage{latexsym}

\newtheorem{theorem}{Theorem}
\newtheorem{lemma}[theorem]{Lemma}
\newtheorem{question}[theorem]{Question}
\newtheorem{definition}[theorem]{Definition}
\newtheorem{conj}[theorem]{Conjecture}

\newcommand {\n}{{\mathbb N}}
\newcommand {\z}{{\mathbb Z}}
\newcommand {\real}{{\mathbb R}}

\begin{document}

\baselineskip=18pt

\begin{center}
{\large\bf Finding integral diagonal pairs in a two dimensional $\mathcal{N}$--set}
\footnotetext{{\em 2000 Mathematics Subject Classification}. Primary 11B75, 11H06, 11P21.}
\end{center}
\begin{center}
{Lev A. Borisov and Renling Jin }
\end{center}
\begin{center}
Abstract
\end{center}
\begin{quote}
\small According to \cite{nathanson} an $n$-dimensional $\mathcal{N}$--set is a
compact subset $A$ of $\real^n$ such that for every $x\in\real^n$ there is
$y\in A$ with $y-x\in\z^n$. We prove that every two dimensional $\mathcal{N}$--set $A$ must
contain distinct points $x,y$ such that $x-y$ is in $\z^2$ and $x-y$ is neither horizontal nor vertical. This answers
a question of P.\ Hegarty and M.\ Nathanson.
\end{quote}

For any sets $A,B$ in an abelian group $A\pm B$ denotes the set $\{a\pm b:a\in A\,\mbox{ and }\,b\in B\}$.
During one of the problem sessions in CANT (Combinatorial and Additive Number Theory Workshop) 2009
M.\ Nathanson asked the following question which was originally raised by P.\ Hegarty:
\begin{question}\label{mainque}
Can we find an $\mathcal{N}$--set $A\subseteq\real^2$, i.e., a compact set $A\subseteq\real^2$ with the property that
$\real^2=A+\z^2$, such that $(A-A)\cap\z^2\subseteq(\z\times\{0\})\cup(\{0\}\times\z)$?
\end{question}

Question \ref{mainque} is motivated by the study of a general inverse problem in order to determine which set
$E\subseteq\z^n$ can be represented by the form of $(A-A)\cap\z^n$ for some $\mathcal{N}$--set $A$.
Notice that $(A-A)\cap\z^n$ contains the origin and is symmetric about the origin.
This inverse problem is completely solved in one dimensional case.  It is shown in \cite{nathanson} that
a finite set $E$ of positive integers is relatively prime if and only if there is an $\mathcal{N}$--set $A\subseteq\real$
such that $E=(A-A)\cap\n$. By the fundamental observation of
geometric group theory (see \cite{nathanson}) if $A$ is an $n$-dimensional $\mathcal{N}$--set, then $(A-A)\cap\z^n$
is a finite set of generators of the group $\z^n$. Clearly, for a one dimensional $\mathcal{N}$--set $A\subseteq\real$,
$(A-A)\cap\z$ is a set of generators if and only if $(A-A)\cap\n$ is relatively prime.
Hence the next logical step is to ask whether a symmetric set of generators of $\z^2$ together with the origin $(0,0)$ can be
represented by the form of $(A-A)\cap\z^2$ for some two dimensional $\mathcal{N}$--set $A$.
For example, it is interesting to ask whether the set $E=\{(0,0),\pm (0,1),\pm(1,0)\}$ can be represented by
$(A-A)\cap\z^2$ for some $\mathcal{N}$--set $A$. The main theorem in this paper will
show that the answer is ``no''. The following is the main theorem.

\begin{theorem}\label{main}
Every $\mathcal{N}$--set $A\subseteq\real^2$ contains $x,y$ such that $x-y\in(\z\smallsetminus\{0\})^2$.
\end{theorem}

\emph{Acknowledgements.} This paper is based on an earlier preprint by the second author, with improvements due to the first author. A somewhat longer argument for the same result has been given independently in \cite{LS}.
The work of the first author was partially supported by the NSF Grant 1003445 and
the work of the second author was partially supported by the NSF Grant RUI 0500671.

\bigskip

We need the following notation.
\begin{definition}
Let $A,B\subseteq\real^2$.
\begin{enumerate}
\item $A$ and $B$ are called
integral lattice congruent provided that $A+z=B$
for some $z\in\z^2$. When we say ``congruent'' in this paper we always mean ``integral lattice congruent.''
\item $A$ and $B$ are called horizontally (vertically) congruent provided that $A+z=B$ for some $z\in\z\times\{0\}$ ($z\in\{0\}\times\z$).
\item $(x,y)\in\real^2\times\real^2$ is called an integral diagonal pair if
$x-y\in\left(\z\smallsetminus\{0\}\right)^2$.
\end{enumerate}
\end{definition}

\noindent {\bf Proof of Theorem \ref{main}}:\quad Assume the contrary and let
$A\subseteq\real^2$ be the $\mathcal{N}$--set which does not contain any integral diagonal pairs.
For each $n\in\n$ let
\[\mathcal{B}_n=\left\{\left[\frac{i}{n},\frac{i+1}{n}\right]
\times\left[\frac{j}{n},\frac{j+1}{n}\right]:i,j\in\z\right\}\] and
\[\mathcal{K}_n(A)=\{B\in\mathcal{B}_n:B\cap A\not=\emptyset\}.\]

Our first step of the proof is to replace the set $A$ by a set $K$, which is the union of finitely many
squares from $\mathcal{B}_n$ for some large $n$. For doing so we transform a continuous problem
to a discrete problem.

\begin{lemma}\label{diag}
There is $N\in\n$ such that for every $n\geqslant N$,
$\bigcup\mathcal{K}_n(A)$ does not contain any integral diagonal pairs.
\end{lemma}

\noindent {\bf Proof of Lemma \ref{diag}}:\quad Since $A$ is compact, we have that $A-A$ is compact.
Since $A$ does not contain integral diagonal pairs, there is $\epsilon>0$ such that
\[\min\{|x-y|:x\in A-A\,\mbox{ and }\,y\in(\z\smallsetminus\{0\})^2\}>\epsilon.\]
Let $N>2\sqrt{2}/\epsilon$. For each $x\in K_n(A)$ where $K_n(A)=\bigcup\mathcal{K}_n(A)$
there is $z\in A$ such that $|x-z|\leqslant\sqrt{2}/n$. Hence
for each $x\in K_n(A)-K_n(A)$ there is a $z\in A-A$ such that
$|x-z|\leqslant 2\sqrt{2}/n$. This implies that for every $n\geqslant N$, for any
$x\in K_n(A)-K_n(A)$ and $y\in(\z\smallsetminus\{0\})^2$ we have that
$|x-y|\geqslant\epsilon-2\sqrt{2}/n>0$. Therefore, $K_n(A)$ does not contain any integral diagonal pairs.
\quad $\Box$

\medskip
We use Lemma \ref{diag} and
fix an $n$ such that $K_n$ contains no integral diagonal pairs.
We will omit the subscript $n$ in
$\mathcal{B}_n$, $\mathcal{K}_n(A)$, etc. Since
$\bigcup\mathcal{K}(A)\supseteq A$, $\bigcup\mathcal{K}(A)$ is compact,
and $\bigcup\mathcal{K}(A)$ does not contain integral diagonal pairs, we could
replace $A$ chosen in the beginning
by $\bigcup\mathcal{K}(A)$ for the counterexample of Theorem \ref{main}.
However, for convenience, we can throw away some unnecessary squares in
$\bigcup\mathcal{K}(A)$. Let
\[k_0=\min\{|\mathcal{K}|:\mathcal{K}\subseteq\mathcal{K}(A)\,\mbox{ and }\,
\bigcup\mathcal{K}\,\mbox{ is an }\mathcal{N}\mbox{--set}\}.\]
\begin{definition}
We call $\mathcal{K}\subseteq\mathcal{B}$ a minimal counterexample of
Theorem \ref{main} if $|\mathcal{K}|=k_0$, $K=\bigcup\mathcal{K}$ is an
$\mathcal{N}$--set, and $K$ does not contain any integral diagonal pair.
\end{definition}

It is easy to see that $k_0=n^2$ because (1) every $B\in\mathcal{B}$ is
congruent to some squares in $\mathcal{K}$ and (2) any two squares
$B,B'\in\mathcal{K}$ are not congruent to each other due to the minimality of $k_0$.
Notice that $A$ may not be covered by a minimal counterexample of Theorem \ref{main}.
However, the assumption that the general counterexample $A$ of Theorem \ref{main} exists implies
the existence of the minimal counterexample $\mathcal{K}\subseteq\mathcal{B}$ of Theorem \ref{main},
we can now forget about $A$ and try to derive a contradiction based on the assumption that
the minimal counterexample $\mathcal{K}\subseteq\mathcal{B}$ of Theorem \ref{main} exists.

We observe that the boundary of the set $K$ consists of a union of horizontal and
vertical segments of length $\frac 1n$ each, which come from the sides of the squares
of $\mathcal{K}$. Each such segment is either horizontally congruent or vertically
congruent to exactly one other segment of the boundary. The kind of congruency
has nothing to do with the direction of the segment itself. The endpoints of each segment
may be congruent to up to four different points on the boundary of $K$ which correspond
to the four adjacent squares on the torus ${\mathbb R}^2/{\mathbb Z}^2$. However,
these points must all lie on one horizontal line or all lie on one vertical line. Conversely,
a collection of $n^2$ small squares, one in each congruency class, whose boundary
satisfies the above properties gives a minimum counterexample. 

Among all possible minimum counterexamples $\mathcal K$ we consider one
with the lowest total length of the boundary of $K$.

Consider the connected components $G_i$ of $\mbox{int}(K)$. To each $G_i$
we associate the area inside the \emph{outer} boundary $\gamma_i$ of $G_i$.
Consider $G_0$ which has the smallest such area among all $G_i$. We claim
that $G_0$ is simply connected. Indeed, otherwise there exists a square $B$,
disjoint from the squares of $G_0$ which lies inside $\gamma_0$. Consider
the region $G_i$ that contains the square $B_1$ congruent to $B$. An integral
shift $G_i+z$ makes it contain $B$, but then $G_i+z$ is a connected region
strictly inside $\gamma_0$. Thus the area inside $\gamma_i$ is smaller
than that inside $\gamma_0$, contradiction.

Now that we know that $G_0$ is simply connected, its boundary is
equal to its outer boundary $\gamma_0$ and is connected. Since each segment
of the boundary is either horizontally or vertically congruent, and the endpoint
of a segment has the same type of congruency, the whole $\gamma_0$
has to have a constant congruency type. Assume that it is horizontally congruent.

Consider a square $B\subseteq G_0$ which lies in the lowest row of $G_0$.
Its bottom $l$ is a part of $\gamma_0$. The square $B_1$ that is adjacent to $B$ at $l$
no longer lies in $G_0$. Moreover, the corresponding congruent square $B_1+ z$
can not lie in $G_0$ by our choice of the row.

Consider now a new set ${\mathcal K}'$ defined by shifting all squares of $G_0$
by $z$ and by keeping the rest of the squares.
We claim that ${\mathcal K}'$ is still a minimum counterexample. Indeed, if we have a
segment on the boundary of $K'=\bigcup {\mathcal K}'$, it could be either a segment
of the boundary of $K$ or a segment of the boundary of $K'$.
Let us first examine its pairs of boundary segments. If they didn't come from $\gamma_0$ or from
segments congruent to $\gamma_0$, then they are unchanged by the switch from
${\mathcal K}$ to ${\mathcal K}'$. If they come from $\gamma_0$ or from
a congruent segment, then the switch may change their relative position. However,
since $G_0$ was shifted horizontally and the relative position was a horizontal congruence,
the new pair of segments might either coincide (thus they are no longer boundary) or
still stay horizontally congruent. Similarly, the only boundary preimages of the corners of the small
squares of the torus that are affected by the change from $\mathcal K$ to
${\mathcal K}'$ are in the congruency class of an endpoint of a segment on $\gamma_0$.
All of these were in a horizontal line in $K$, so they stay in a horizontal line in $K'$.
If in fact they all collapse to a point, then we no longer have a point on the boundary,
but certainly, we do not introduce any integral diagonal pairs.

It remains to observe that $l$ is no longer a boundary of $K'$, so the length
of the boundary of $K'$ is less than that of $K$, which contradicts the
minimality in the choice of $\mathcal K$.
\quad $\Box$

\bigskip

The proof of Theorem \ref{main} can be easily modified to prove the
following generalization.

\begin{theorem}\label{main2}
If $A\subseteq\real^2$ is an $\mathcal{N}$--set and $l_1,l_2$ are two one dimensional linear subspaces
of $\real^2$, then $(A-A)\cap\z^2\not\subseteq l_1\cup l_2$
\end{theorem}

\medskip

Although we believe that the multi-dimensional version of Theorem \ref{main} should be
true, we cannot find a proof. Thus we would like to make the following conjecture.

\begin{conj}
If $A\subseteq\real^3$ is an $\mathcal{N}$--set and $p_1,p_2,p_3$ are three two dimensional linear
subspaces of $\real^3$, then $(A-A)\cap\z^3\not\subseteq p_1\cup p_2\cup p_3$.
\end{conj}

\bigskip

\noindent {\em Department of Mathematics, Rutgers University, Piscataway, NJ 08854, USA}

\medskip

\noindent {\em Department of Mathematics, College of Charleston, Charleston, SC 29424, USA}

\end{document}